\documentclass[a4paper]{article}
\usepackage{amssymb}
\usepackage{amsmath}
\usepackage{amsfonts}
\usepackage{graphicx}
\usepackage{pstricks}
\begin{document}
\title{Homotopy type of the complement of an immersion and
classification of embeddings of tori}
\author{M. Cencelj, D. Repov\v s and M. Skopenkov
\thanks{Cencelj and Repov\v s were supported in part by the Slovenian
Research Agency Program No.~101-509. Skopenkov was supported in
part by President of the Russian Federation grant NSh-4578.2006.1,
Agency for Education and Science grant RNP-2.1.1.7988,
Russian Foundation for Basic Research grants
No.~02-01-00014, 05-01-00993, 06-01-72551,
07-01-00648, and INTAS grant No.~06-1000014-6277.}}
\date{}
\maketitle

This paper is devoted to the classification of embeddings of
higher dimensional manifolds. This subject was actively studied in
the sixties  [1] ,  [2]  and there has been a
renewed interest in it in  recent years  [3] -- [5].
Investigation of this problem started with consideration of
knots $S^q\to S^m$, for which  an
explicit classification in some dimensions, and a complete
{\it rational} classification in codimension $\ge 3$  was obtained:

\smallskip
\noindent{\bf Theorem 1.}  [1]  {\it
Assume that $q+2 < m < \frac{3}{2}q+2$. Then up to isotopy the
set of smooth embeddings $S^q\to S^m$ is infinite if and only if
$q+1$ is divisible by $4$.}
\smallskip

We study the case of embeddings $S^p\times S^q\to S^m$.
We call such embeddings {\it knotted tori}. Classical
special case of knotted tori are links (see Figure 1a). The
investigation of knotted tori is a natural next step after knots and links
because of the handle decomposition of an arbitrary
manifold.

\begin{figure}[ht]
\center
\includegraphics{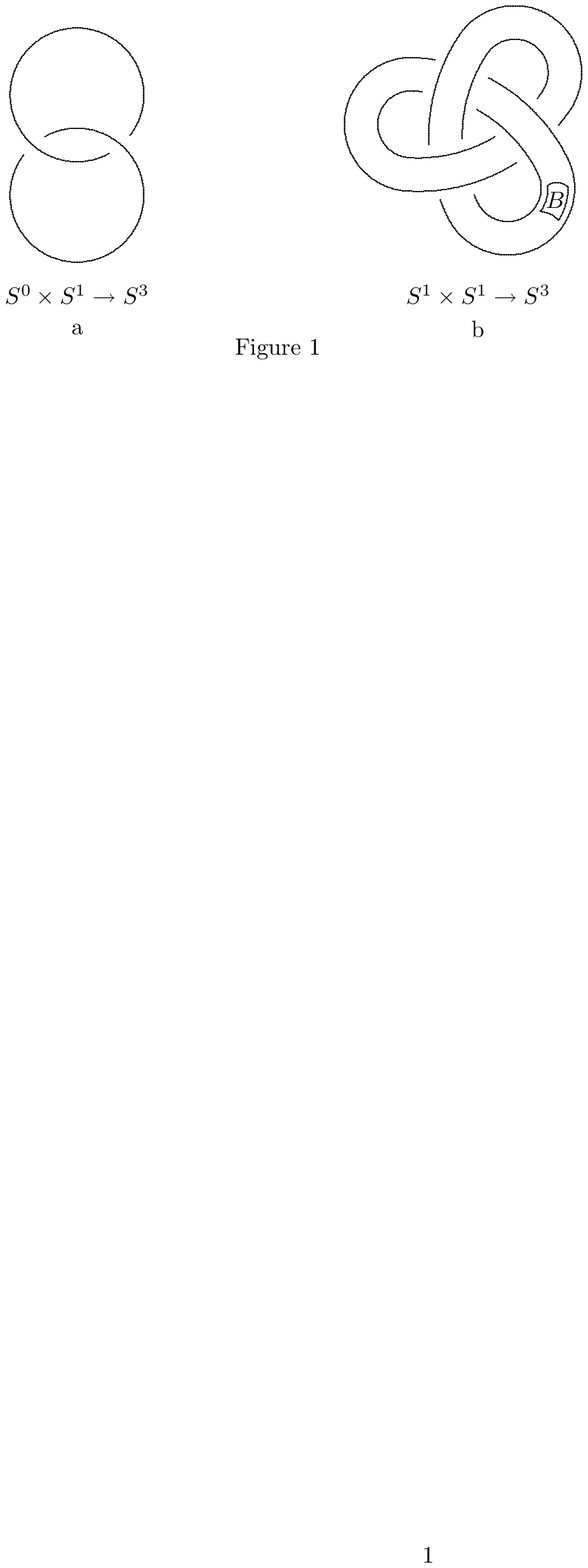}
\end{figure} 

The set of knotted tori in the the space of sufficiently high dimension, namely in  the {\it metastable} range $m\ge p+3q/2+2$, $p\le q$,
was explicitely described in  [4].
The metastable range is a natural limit
for the
classical methods of embedding theory.
The aim of this note is to present an approach which allows for  results
in lower dimension:

\smallskip
\noindent{\bf Theorem 2.} {\it Assume that $p+4q/3+2 < m <
p+3q/2+2$ and $m > 2p+q+2$. Then the set of smooth
embeddings $S^p\times S^q\to S^m$ up to smooth isotopy is infinite if and
only if either $q+1$ or $p+q+1$ is divisible by $4$.}
\smallskip

Our approach to classification is based on investigation of immersions and the homotopy type of their complements.

Let us give some necessary definitions. Fix points $u\in S^p$ and $v\in S^q$
and a ball $B=B^{p+q}\subset S^p\times S^q$ 
which does not intersect the meridian $S^p\times v$ 
and the parallel $u\times S^q$ of the torus (see Figure~1b).
A piecewise smooth map $F:S^p\times S^q\to S^m$ is called
an {\it almost
embedding}, if $F$ is a smooth embedding outside the
ball $B$, and $FB\cap F(S^p\times S^q\setminus B)=\varnothing$.
An {\it almost isotopy} is defined analogously, only $B$ is replaced by $B\times I$.
Denote by $\mathcal M$ the set of all almost embeddings $F:S^p\times
S^q\to S^m$ up to almost isotopy.
As it was shown in  [5]  the set $\mathcal M$ has a natural group structure.

The description of the group $\mathcal M$ is a much simpler problem than
the classification of embeddings and it is done by classical methods.
Thus Theorem~2 is reduced to the following problem: to
{\it determine, when a given almost embedding is almost isotopic to an embedding}.
Our main lemma asserts that the complete obstruction to this lies in the finite group:

\smallskip
\noindent{\bf Lemma 3.} {\it Let  $p+4q/3+2 < m <
p+3q/2$ and $m > 2p+q+2$. Then there is a homomorphism $\beta: \mathcal M\to
\mathcal G$ into a finite group $\mathcal G$,
with the following property: if $\beta(F)=0$ then $F$ is almost isotopic to an embedding
(smooth outside~$B$).}
\smallskip
 
We give a sketch the proof of Theorem 2 in the case when $q+1$ is divisible by 4
under assumption that Lemma~3 has been  proved. It can be shown that the group
$\mathcal G$ is infinite in this case (and finitely generated).
Then the group $|\mathcal G|\cdot \mathcal M$ is also infinite.
For each $F\in |\mathcal G|\cdot \mathcal M$ we have $\beta(F)=0$,
thus each such $F$ 
 is almost isotopic to an embedding (smooth outside the ball $B$).
By smoothing these embeddings we get infinitely many distinct smooth embeddings $S^p\times S^q\to S^m$.

\smallskip
{\bf Proof} of Lemma~3 is based on the following 3-step construction (see Figure~2):

\begin{figure}[ht]
\center
\includegraphics{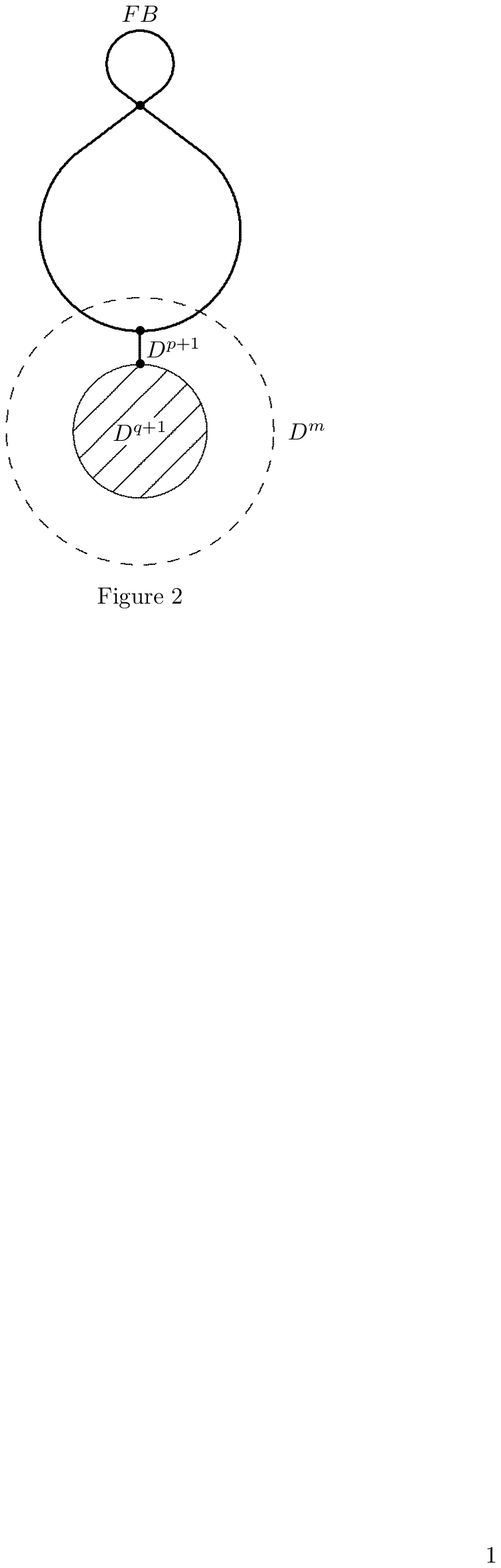}
\end{figure} 

{\it Step 1. Construction of a web $D^{p+1}$.} Let  $F:S^p\times S^q\to S^m$  be an almost embedding.
Glue a disc $D^{p+1}\subset S^m$ to the meridian $F(S^p\times v)$, whose interior is missing the
rest of the  torus $F(S^p\times S^q)$. 
This is possible by general position, because $m> 2p+q+2$. We call such a disc $D^{p+1}$ is a {\it  web}.

\smallskip
{\it Step 2. Construction of a web $D^{q+1}$.} Using general position we cannot construct the second web $D^{q+1}\subset S^m$ 
glued along the parallel $F(u\times S^q)$. So let us assume initially that
such a web exists.

\smallskip
{\it Step 3. Conical Construction.} Remove from the sphere $S^m$ 
a neighborhood of the union of the webs. By general position generically the webs have a unique common point, therefore as a result we obtain an $m$-dimensional ball $D^m$. 
We may assume that $D^m\cap F(S^p\times S^q)=FB$.   
By replacing $FB$ by the cone over $F\partial B$ inside $D^m$, we transform $F$ to the required embedding. 
\smallskip

The outline of the rest of the proof of Lemma~3 is as follows. The desired homomorphism $\beta: \mathcal M\to
\mathcal G$ is the obstruction to existence of the web $D^{q+1}\subset S^m$. This obstruction lies in the group $\pi_q(D^m-\mathop{Im}F,\partial D^m-\mathop{Im}F)$, where $D^m$ is the complement of
a neighborhood of the web $D^{p+1}$. This group is evaluated by the methods of [3; \S4], and in our case  it is
finite. 

\bigskip

\noindent[1] N. Habegger and U. Kaiser,
Topology 37:1 (1998), p. 75--94.

\noindent[2] A. Haefliger,
Ann. Math. 83:3 (1966), p. 402--436.

\noindent[3] A.~Haefliger,
Comment. Math. Helv. 41 (1966-67), p. 51--72.

\noindent[4] A. Skopenkov,
Comment. Math. Helv. 77 (2002), p. 78--124.

\noindent[5] A. Skopenkov, Classification of embeddings below the
metastable dimension, preprint, http://arxiv.org/ math.GT/0607422

\end{document}